\documentclass{article}
\usepackage{enumitem}

\usepackage{amsmath}
\usepackage{amsthm}
\usepackage{graphicx}
\newtheorem{example}{Example}
\newtheorem{proposition}{Proposition}
 \usepackage[main, final]{neurips_2025}

\usepackage[utf8]{inputenc} 
\usepackage[T1]{fontenc}    
\usepackage{hyperref}       
\usepackage{url}            
\usepackage{booktabs}       
\usepackage{amsfonts}       
\usepackage{nicefrac}       
\usepackage{microtype}      
\usepackage{xcolor}         

\title{GPU-based Split algorithm for Large-Scale CVRPSD}

\author{%
  Jingyi Zhao\textsuperscript{1}, Linxin Yang\textsuperscript{1,2}, Haohua Zhang\textsuperscript{4}, Tian Ding\textsuperscript{1,3,4}\\
  1. Shenzhen Research Institute of Big Data, Shenzhen, China\\
  2. The Chinese University of Hongkong, Shenzhen, China\\
  3. Shenzhen International Center for Industrial and Applied Mathematics, Shenzhen, China\\
  4. AutoKernel, Shenzhen, China
}

\begin{document}

\maketitle

\begin{abstract}
Dynamic programming (DP) is a cornerstone of combinatorial optimization, yet its inherently sequential structure has long limited its scalability in scenario-based stochastic programming (SP). This paper introduces a GPU-accelerated framework that reformulates a broad class of forward DP recursions as batched min--plus matrix--vector products over layered DAGs, collapsing actions into masked state-to-state transitions that map seamlessly to GPU kernels. Using this reformulation, our approach takes advantage of massive parallelism across both scenarios and transitions, enabling the simultaneous evaluation of \emph{over one million uncertainty realizations} in a single GPU pass -- a scale far beyond the reach of existing methods. We instantiate the framework in two canonical applications: the capacitated vehicle routing problem with stochastic demand and a dynamic stochastic inventory routing problem. In both cases, DP subroutines traditionally considered sequential are redesigned to harness two- or three-dimensional GPU parallelism. Experiments demonstrate near-linear scaling in the number of scenarios and yield one to three orders of magnitude speedups over multithreaded CPU baselines, resulting in tighter SAA estimates and significantly stronger first-stage decisions under fixed time budgets. Beyond these applications, our work establishes a general-purpose recipe for transforming classical DP routines into high-throughput GPU primitives, substantially expanding the computational frontier of stochastic discrete optimization to the million-scenario scale.
\end{abstract}

\section{Introduction}
\paragraph{Background.}

The Capacitated Vehicle Routing Problem (CVRP, \citep{toth2002models,laporte2009fifty,pecin2017improved,accorsi2021fast,vidal2022hybrid}) is a classical NP-hard problem with broad industrial relevance in logistics, last-mile delivery, and supply chains.
The task is to design cost-minimizing routes for a fleet of identical, capacity-limited vehicles that start and end at a central depot, ensuring that all customers are served exactly once while respecting vehicle capacity constraints and satisfying customer demand.
Its stochastic extension, the CVRPSD (\citep{laporte2002integer,gauvin2014branch,florio2020new,fukasawa2023complexity,ota2025hardness}), models customer demands as random variables within a two-stage program: routes are chosen before demand revelation, followed by recourse actions such as restocking or rerouting. Despite extensive research (\citep{rei2007local,gounaris2013robust,mendoza2011constructive}), solving CVRPSD remains challenging in data-driven settings with correlated or complex demand distributions. State-of-the-art set-partitioning approaches often rely on restrictive independence assumptions to stay tractable.
Scenario-based formulations offer greater modeling flexibility, but evaluating large scenario sets is computationally prohibitive on CPUs.

Stochastic programming naturally lends itself to parallelism: given a first-stage decision, the evaluation of second-stage costs across different scenarios is inherently independent.
Thus, we propose a GPU-accelerated framework that reformulates the split algorithm into a vectorized, GPU-parallelizable form, enabling simultaneous evaluation across scenarios and transitions. This design supports millions of scenarios within practical runtimes, substantially reducing estimation bias in sample-average approximation.
While researchers in both operations research and machine learning have increasingly explored GPU-based parallelism to accelerate optimization algorithms (e.g., GPU-ADMM~\citep{schubiger2019gpu},CuPDLP~\citep{lu2023cupdlp},  GPU-QP~\citep{bishop2024relu}, GPU-IPM~\citep{liu2024gpu}), to the best of our knowledge, this work is the first to introduce an ultra-large-scale scenario-based solver for stochastic combinatorial optimization, capable of efficiently handling \emph{over 1,000,000 of uncertainty realizations}. This represents a significant advancement in the computational frontier of scenario-based stochastic programming, which has traditionally been constrained by the computational cost of evaluating massive scenario sets.


\paragraph{Contribution.}
Our framework evaluates a first-stage routing decision by computing the second-stage recourse costs for all demand scenarios simultaneously on the GPU. This parallel evaluation enables rapid and reliable comparison of candidate solutions across ultra-large-scale scenario sets, providing an almost unbiased estimate of the expected cost. Consistent with the theoretical insight that larger scenario sets reduce estimation bias, our experimental results confirm that first-stage solutions derived from richer scenario sets achieve superior performance on out-of-sample tests.

To the best of our knowledge, this is the first approach capable of evaluating over 1,000,000 of demand scenarios within practical runtimes for first-stage decision making. More broadly, it represents one of the earliest systematic applications of GPU-based parallelization to stochastic combinatorial optimization, moving beyond the conventional reliance on CPU-based multi-threading. 
Our results further demonstrate that GPU parallelization enables much finer-grained computation than scenario-level simulation, achieving speedups of up to $65\times$ over CPU-based multi-threading and thereby unlocking unprecedented scalability.

Our parallel evaluation strategy integrates seamlessly into the Hybrid Genetic Search (HGS) framework originally proposed by~\citep{vidal2012hybrid}, a widely used metaheuristic that has achieved state-of-the-art performance across numerous VRP variants~\citep{vidal2013hybrid,vidal2014unified, vidal2015hybrid, goel2021team}. Importantly, the GPU-based parallel sub-problem solution module is generic: it can be embedded not only into metaheuristics but also into exact and approximate algorithms that require large-scale scenario evaluation. Beyond the CVRPSD, our parallel simulation approach is readily extensible to dynamic programming and rollout-based methods under massive scenario sets, offering a scalable tool for a broad class of stochastic optimization problems. 

Viewed from the perspective of stochastic programming, our framework bridges a key methodological gap by substantially alleviating the long-standing scalability barrier in data-driven combinatorial optimization, thereby opening the door to previously intractable real-world applications.

\section{Problem Setting}
The stochastic programming community has extensively studied scenario-based formulations, where uncertainty is modeled by a finite set of realizations. 
However, scenario-based evaluation quickly becomes computationally prohibitive on CPUs, where even tens of thousands of scenarios can overwhelm multi-threaded implementations.
In our work, we adopt such a \emph{scenario-based modeling} framework, in which customer demands are represented by sampled realizations. This approach naturally accommodates correlated demand structures and supports data-driven modeling when historical records are available. This scenario-based approach falls naturally into the two-stage stochastic programming paradigm, whose general form is: $\min_{x \in X} f_1(x) + \mathbb{E}_{\xi} \left[ f_2(x, \xi) \right],$ where $x$ denotes the first-stage routing decisions and the corresponding cost $f_1(x)$, $\xi$ is a random vector representing a realization of customer demands, and $f_2(x, \xi)$ denotes the second-stage recourse cost under each scenario $\xi$.

In our two-stage stochastic optimization problem, the first stage determines the visiting sequence of customers, commonly referred to as a \emph{giant tour} in the context of genetic algorithms~\citep{vidal2012hybrid}. This representation encodes a solution as a permutation of customers, from which feasible vehicle routes can be recovered via a split operator under fixed vehicle capacity and scenario-dependent customer demands. We assume \emph{full demand revelation prior to the second stage}, enabling the plan to adapt to the realized scenario. Once demands are known, the giant tour is \emph{split into feasible routes} such that the demand on each route does not exceed vehicle capacity. Thus, given a fixed giant tour (i.e., the visiting sequence), the second-stage evaluation is computationally simple: it only requires splitting the tour according to realized demands. \textbf{The objective is to determine a first-stage giant tour that minimizes the expected total travel cost across all scenarios.}

\begin{figure}[htbp]
    \centering
    \includegraphics[width=0.85\linewidth]{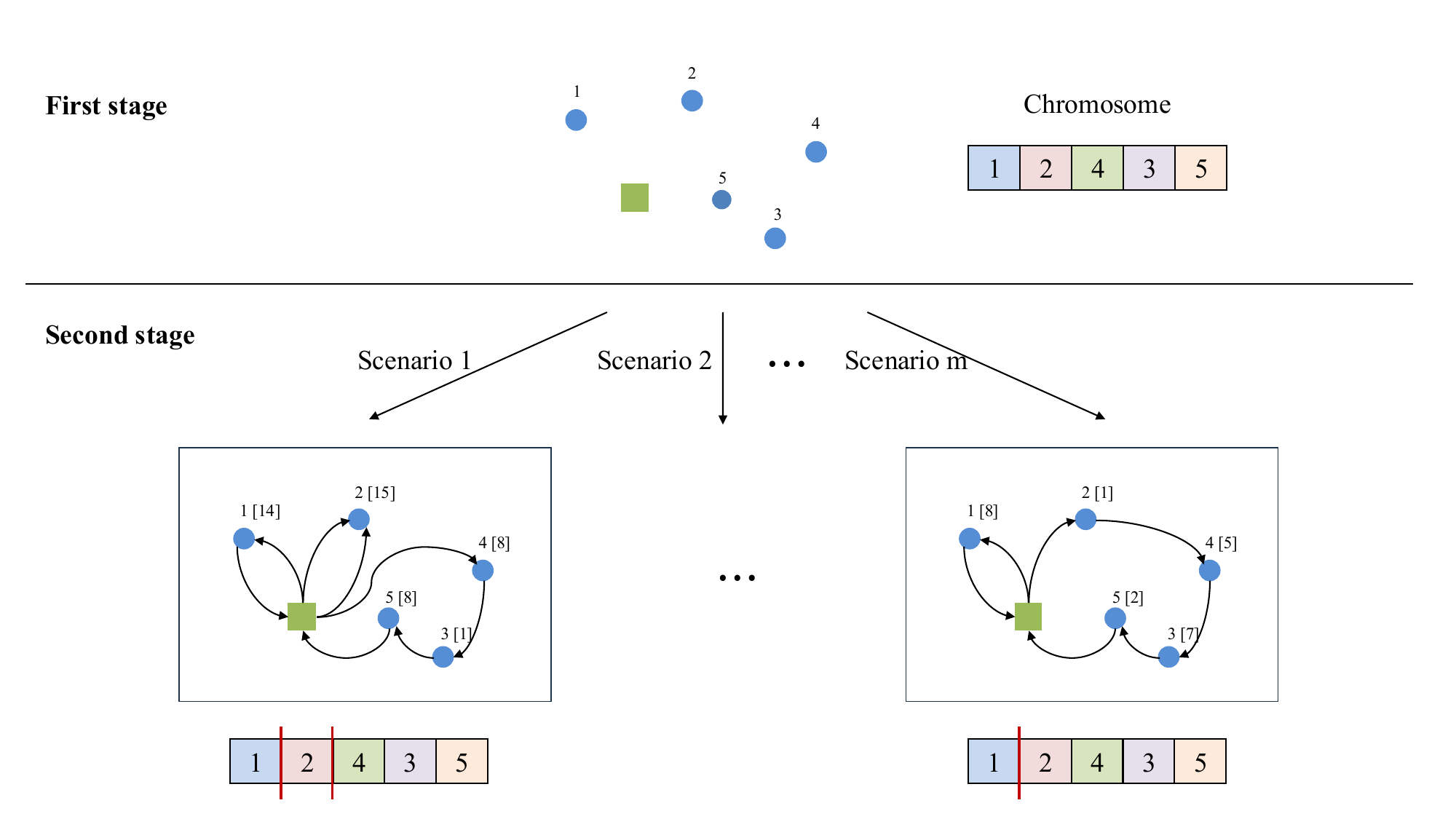}
    \caption{\small{Example of splitting a giant tour into feasible routes under a given demand scenario.}}
    \label{fig:split}
\end{figure}
\begin{example}
Suppose the first-stage giant tour is $(1,2,4,3,5)$ with vehicle capacity $Q=17$, 
as illustrated in Figure~\ref{fig:split}. 

\begin{enumerate}[label=(\roman*)]
    \item In Scenario~1, the realized demands are $[14,15,8,1,8]$. 
    The second-stage routing may then split the tour into three feasible routes: 
    $(1)$, $(2)$, and $(4,3,5)$. 

    \item In Scenario~$m$, the realized demands are $[8,1,5,7,2]$ with three vehicles available. 
    The second-stage routing may then split the tour into two feasible routes: 
    $(1)$ and $(2,4,3,5)$. 
\end{enumerate}

\end{example}
 It is well understood that training robust first-stage policies benefits from incorporating a large number of demand scenarios to capture the full breadth of uncertainty (\cite{shapiro2003monte}, see details in the Appendix).
 Theoretical results show that a sufficiently large scenario set is critical for achieving accurate and stable estimations. 
However, in practice, evaluating even tens of thousands of scenarios can be prohibitively expensive, especially for combinatorial problems such as the CVRPSD, since each scenario requires solving a non-trivial routing evaluation (typically through a dynamic programming algorithm with time complexity $O(nP)$, where $n$ is the number of customers and $P$ is the number of possible transitions) rather than a simple function call.



\section{GPU-Accelerated Dynamic Programming-based Split Algorithm}
\label{sec:GPU}
\paragraph{Dynamic Programming-based Second-stage Split Algorithm.}
Let a first-stage decision \emph{giant tour} be denoted by $\sigma = [\sigma_0, \sigma_1, \ldots, \sigma_n, \sigma_{n+1}]$, where $\sigma_0 = 0$ is the departure depot, $\sigma_{n+1} = n+1$ is the return depot, and $\sigma_1, \ldots, \sigma_n$ are the customers in the first-stage visiting sequence.

In the second-stage, for a given demand scenario $\omega \in \Omega$ with realized demand vector $\xi^\omega = (q_{\sigma_1}^\omega, \ldots, q_{\sigma_n}^\omega)$, the task is to \emph{split} the giant tour into a set of feasible vehicle routes that respect both the vehicle capacity $Q$ and, if applicable, maximum route length or duration constraints, while minimizing the total travel cost.

We solve this splitting problem using a dynamic programming (DP) formulation.  
Let $f^\omega(i)$ denote the minimum total cost of serving customers $[\sigma_1, \ldots, \sigma_i]$ starting from the depot and returning to it under scenario $\omega$.  
The base case is $f^\omega(0) = 0$.  
For $i \geq 1$, $f^\omega(i)$ is computed recursively from preceding states:
\begin{equation}
\label{eq:dp-split-scenario}
f^\omega(i) = \min_{0 \leq p \leq i-1} \left\{ f^\omega(p) + c_{0,\sigma_{p+1}} + \sum_{k=p+1}^{i-1} c_{\sigma_k, \sigma_{k+1}} + c_{\sigma_i, n+1} \;\middle|\; \sum_{k=p+1}^{i} q_{\sigma_k}^\omega \leq Q \right\},
\end{equation}
where $c_{a,b}$ denotes the travel cost from node $a$ to node $b$.
Here, $p$ represents the index of the last customer in the previous route, so that $\sigma_{p+1}$ is the first customer in the current route.  
The capacity constraint in~\eqref{eq:dp-split-scenario} ensures that the total demand from $\sigma_{p+1}$ to $\sigma_i$ in scenario $\omega$ does not exceed $Q$.  
The recursion effectively enumerates all feasible last segments of the solution, selecting the one that yields the minimal total cost.

\paragraph{GPU-Parallelized Vectorized Split Algorithm.}
\begin{figure}[htbp] 
\hspace{-1.5cm}
\centering 
\includegraphics[width=0.85\linewidth]{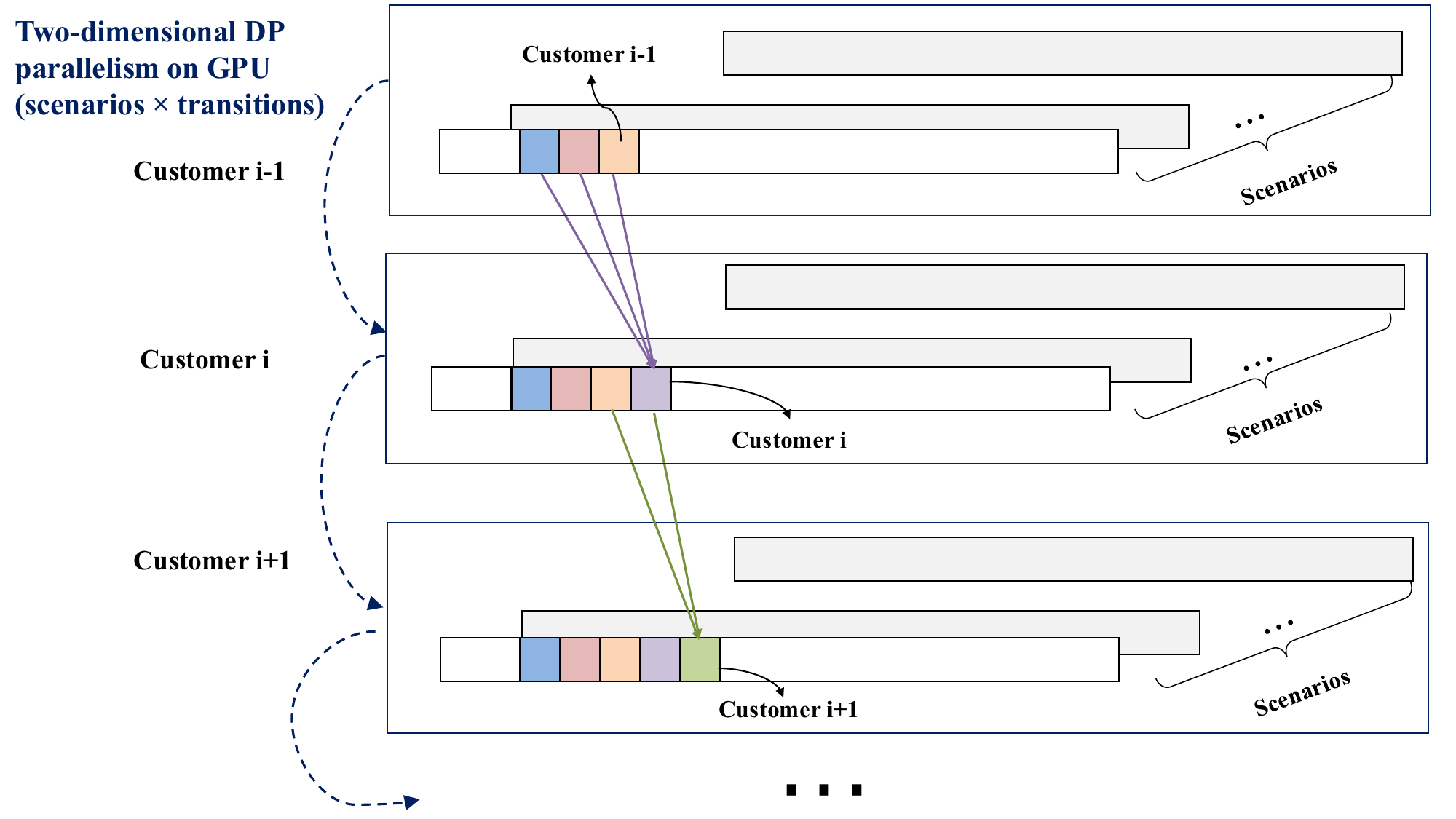} \caption{\small{Example of Vectorized Split Algorithm.} }\label{fig:matrix} 
\end{figure}
While the vanilla split algorithm in~\eqref{eq:dp-split-scenario} is straightforward, it requires
checking the capacity constraint for every candidate split point $p$ at runtime.  
On a GPU, such per-transition feasibility checks can cause thread divergence and redundant memory access, 
severely limiting throughput when evaluating a large number of scenarios in parallel.

To address this, we introduce a \emph{vectorized split algorithm} that precomputes, for each scenario $\omega \in \Omega$ 
and each position $i$, the earliest feasible split point satisfying the capacity constraint:
\begin{equation}
\label{eq:mask-def}
\text{mask}^\omega(i) \;=\; \min \left\{ p \,\middle|\, 0 \leq p < i,\; \sum_{k=p+1}^i q_{\sigma_k}^\omega \leq Q \right\}.
\end{equation}
Here, $\text{mask}^\omega(i)$ denotes the smallest index $p$ such that the segment 
$[\sigma_{p+1}, \ldots, \sigma_i]$ is feasible under scenario $\omega$.  
This precomputation can be carried out in parallel across all scenarios and all customer positions 
using prefix-sum operations on the demand vectors $\xi^\omega$.

With $\text{mask}^\omega(i)$ available, the DP recursion becomes:
\begin{equation}
\label{eq:dp-split-mask}
f^\omega(i) = \min_{\text{mask}^\omega(i) \leq p \leq i-1} 
\left\{ f^\omega(p) + c_{0,\sigma_{p+1}} + \sum_{k=p+1}^{i-1} c_{\sigma_k, \sigma_{k+1}} + c_{\sigma_i, n+1} \right\}.
\end{equation}
The capacity feasibility check is now implicit in the lower bound $\text{mask}^\omega(i)$, eliminating
conditional branching inside the DP loop.

The vectorized formulation naturally enables parallelization along two dimensions, as shown in Figure~\ref{fig:matrix}:
\begin{enumerate}
    \item \emph{Scenario-level parallelism:} each demand scenario $\omega$ is processed independently, 
    making it straightforward to evaluate thousands of scenarios simultaneously on different GPU threads.  
    \item \emph{Transition-level parallelism:} for each state $i$, the minimization over candidate split points 
    $p \in [\text{mask}^\omega(i), i-1]$ can itself be parallelized, since all candidate costs are independent 
    and can be reduced to a minimum using parallel reduction.
\end{enumerate}

These two levels of parallelism substantially improve GPU efficiency:
(1) all threads follow the same execution path, reducing warp divergence, and 
(2) memory accesses are more regular since $\text{mask}^\omega(i)$ values are stored contiguously, 
allowing coalesced reads.  
In practice, this reduces per-scenario runtime significantly, enabling the evaluation of tens or 
hundreds of thousands of scenarios within practical time limits.

\section{Experiments}

    We conduct preliminary experiments to evaluate the efficiency of our GPU-accelerated split algorithm against CPU baselines. All implementations were written in C++/CUDA, and experiments were run on a machine with an AMD Ryzen 7 9700X CPU (8 cores) and an NVIDIA RTX 2080Ti GPU with 11 GB memory. The CPU baselines include (i) a single-threaded implementation, and (ii) a multi-threaded implementation with 8 threads. The GPU implementation exploits the two-dimensional parallelism described in Section~\ref{sec:GPU}. 

\subsection{Scalability with the Number of Scenarios}

\begin{figure}[htbp]
    \centering
    \includegraphics[width=0.55\linewidth]{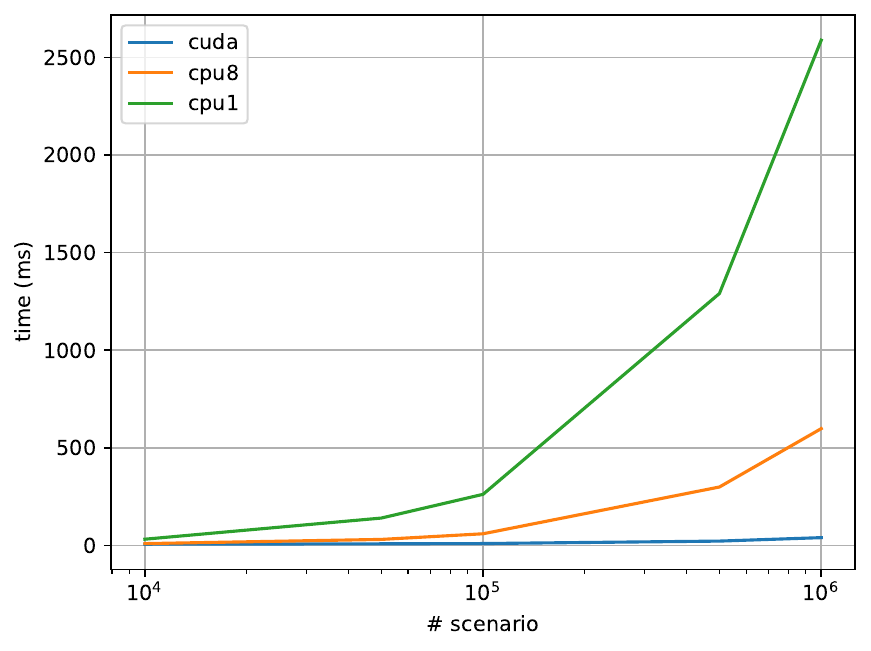}
    \caption{\small{Runtime comparison for processing a splits evaluation under $10^4$–$10^6$ demand scenarios. GPU parallelization scales almost linearly with the number of scenarios, while CPU runtimes grow super-linearly and quickly become prohibitive.}}
    \label{fig:runtime-vs-scenario}
\end{figure}

We first compare the runtime required to split a single giant tour under varying numbers of considered scenarios, ranging from $10^4$ to $10^6$. 
The GPU-based split algorithm is compared against single thread and 8-threads CPU baseline.
Figure~\ref{fig:runtime-vs-scenario} presents the average runtime (in milliseconds) computed over 10 independent runs.

\paragraph{Results.} 
On the CPU single-thread baseline, the runtime grows rapidly and exceeds several minutes when $10^6$ scenarios are evaluated.  
The multi-threaded CPU implementation achieves moderate improvements, but scaling saturates due to synchronization and memory bandwidth bottlenecks.  
In contrast, the GPU implementation achieves near-linear scaling with the number of scenarios, maintaining runtimes within a few seconds even at $10^6$ scenarios.  

Overall, the GPU implementation provides \emph{two to three orders of magnitude speedups} over CPU single-thread and more than one order of magnitude over CPU multi-threading. This validates our claim that GPU parallelism fundamentally shifts the computational frontier for scenario-based stochastic evaluation.

\subsection{Impact of Training Scenario Set Size on Decision Quality}
\begin{figure}[htbp]
    \centering
    \makebox[\textwidth][c]{   
        \begin{minipage}{0.48\textwidth}
            \centering
            \includegraphics[width=\linewidth]{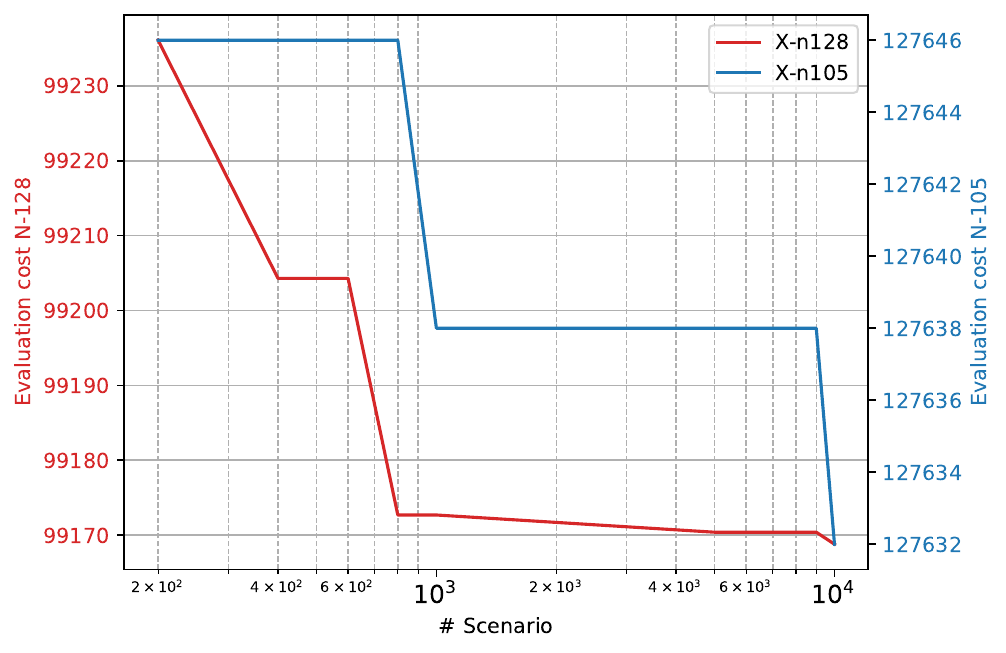}
    \caption{\small{Out-of-sample performance of first-stage solutions obtained with varying observed scenario settings. Larger  evaluation set yield more robust and lower-cost solutions.}}
    \label{fig:trainset-vs-cost}
        \end{minipage}
        \hspace{0.05\textwidth}
        \begin{minipage}{0.46\textwidth}
            \centering
            \includegraphics[width=\linewidth]{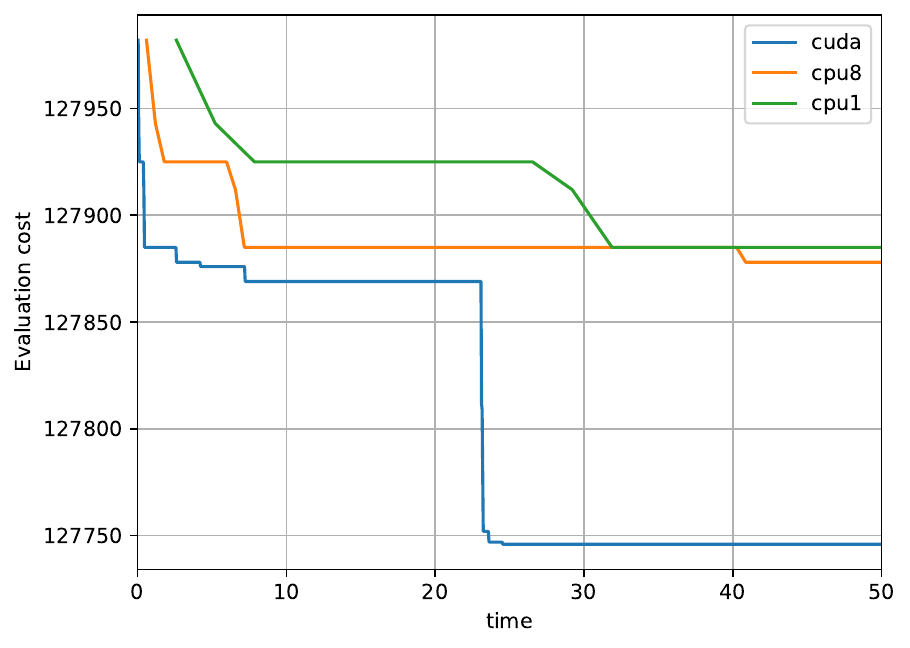}
    \caption{\small{Quality of the best solution obtained at each time under a fixed time budget. GPU consistently achieves better decisions due to faster evaluation and thus larger effective search effort.}}
    \label{fig:time-vs-quality}
        \end{minipage}
    }
\end{figure}

We next examine how the number of evaluated scenarios influences the quality of the first-stage decision. Specifically, we solve the problem under different scenario counts, ranging from 1 to $10^4$ (i.e., $1, 100, 1{,}000$).
For each scenario setting, the obtained first-stage solution is evaluated on a fixed large out-of-sample test set of $10^6$ scenarios. 
Figure~\ref{fig:trainset-vs-cost} reports the out-of-sample cost achieved by the best observed solution on two CVRPSD instances: x-n128 and x-n105.

\paragraph{Results.} 
When observing only few scenarios (e.g., a single scenario), the resulting first-stage solution is severely biased and performs poorly.
Increasing the scale of available scenarios consistently improves robustness, with significant gains observed throughout Figure~\ref{fig:trainset-vs-cost}.
Such performance confirms the theoretical insight that larger sample sizes reduce estimation bias in sample-average approximation.
These results demonstrate that our GPU-based framework, by enabling the evaluation of tens or hundreds of thousands of scenarios within practical runtimes, leads to significantly more reliable first-stage solutions compared to CPU-based methods that are restricted to only a few thousand scenarios.

\subsection{Decision Quality under Fixed Time Budgets}

Finally, we compare the decision quality obtained under identical wall-clock time limits across the three implementations: CPU single-thread, CPU multi-thread, and GPU.  
Each method is given a fixed runtime budget , during which the split algorithm splits giant tours to obtain first-stage solutions.
For fairness, all approaches are evaluated on the same problem instance with $10^4$ available scenarios.  
Figure~\ref{fig:time-vs-quality} reports the best penalized cost obtained within the allowed runtime.

\paragraph{Results.} 
With small time budgets, all methods return feasible but suboptimal solutions, yet GPU already provides a noticeable advantage.  
As the time limit increases, the quality gap widens: GPU produces solutions that are consistently closer to the true optimum, while CPU single-thread stagnates and CPU multi-thread improves only modestly.  
This matches intuition: faster scenario evaluation allows the GPU to explore many more candidate first-stage tours within the same runtime, thereby improving the probability of discovering high-quality solutions.  

These results confirm that beyond scalability, GPU acceleration directly translates into superior decision quality under realistic time constraints, making it particularly valuable in operational settings where decisions must be made quickly.

\section{Conclusion}
We proposed a GPU-accelerated framework for scenario-based stochastic combinatorial optimization, exemplified by the CVRPSD. By reformulating the split algorithm into a vectorized GPU-compatible form, our approach achieves orders-of-magnitude speedups over CPU baselines while enabling evaluation on ultra-large scenario sets. Experiments show that GPU parallelism not only scales efficiently but also yields more robust first-stage solutions and superior decision quality under fixed time budgets. These results demonstrate that GPU acceleration fundamentally advances the scalability and effectiveness of scenario-based stochastic programming, with broad applicability to logistics and other data-driven optimization domains.
\section{Acknowledge}
The work of Jingyi Zhao is supported by the Shenzhen Research Institute of Big Data under grant number J00220250002.
The work of Tian Ding is supported by Hetao Shenzhen-Hong Kong Science and Technology Innovation Cooperation Zone Project (No. HZQSWS-KCCYB-2024016) and Natural Science Foundation of China (Grant No. 12401409).
\small
\bibliographystyle{plain}
\bibliography{sample}

{


\appendix

\section{Monte Carlo Method Proposition}

\emph{Empirical Risk Minimization} (ERM) method is analogous to the Monte Carlo method for estimating a population mean via sample averages and fits naturally within the ERM framework for stochastic programming.  
Formally, we distinguish between:

\begin{itemize}
    \item \emph{True problem}:
    \begin{equation*}
        \small
        (P) \quad z^* = \min_{x \in \mathbb{X}} \mathbb{E}[f(x, \tilde{\xi})],
    \end{equation*}
    
    \item \emph{Sample-average problem} with \( m \) scenarios:
    \begin{equation*}
        \small
        (P_m) \quad z_m^* = \min_{x \in \mathbb{X}} \frac{1}{m} \sum_{i=1}^{m} f(x, \tilde{\xi}^i),
    \end{equation*}
\end{itemize}

As the sample size $m$ increases, the optimal solution $x_m^*$ of the sample problem converges to the true optimal solution $x^*$, and the optimal value $z_m^*$ approaches $z^*$. The following result summarizes the fundamental properties of the ERM method under mild regularity conditions (cf.~\citep{shapiro2003monte}).

\begin{proposition}
Let $\{ \tilde{\xi}^1, \dots, \tilde{\xi}^m \}$ be i.i.d. samples of $\tilde{\xi}$. Denote by
\[
    (P) \quad z^* = \min_{x \in \mathbb{X}} \mathbb{E}[f(x, \tilde{\xi})], 
    \qquad
    (P_m) \quad z_m^* = \min_{x \in \mathbb{X}} \frac{1}{m} \sum_{i=1}^m f(x, \tilde{\xi}^i),
\]
the true and sample average problems, respectively, with optimal solutions $x^*$ and $x_m^*$. Then:
\begin{align}
    &\text{(Bias)} && \mathbb{E}[f(x_m^*, \tilde{\xi})] \leq z^*, \label{eq:bias} \\
    &\text{(Consistency)} && \mathbb{E}[f(x_m^*, \tilde{\xi})] \xrightarrow{\text{a.s.}} z^* \quad \text{as } m \to \infty, \label{eq:consistency} \\
    &\text{(Probabilistic Convergence)} && \lim_{m \to \infty} \Pr\left\{ \mathbb{E}[f(x_m^*, \tilde{\xi})] \leq z^* + \tilde{\epsilon}_m \right\} \geq 1 - \alpha, 
    \quad \tilde{\epsilon}_m \downarrow 0, \label{eq:probconv} \\
    &\text{(Rate of Convergence)} && \sqrt{m}\,(z_m^* - z^*) \xrightarrow{d} \mathcal{N}\big(0, \sigma^2(x^*)\big). \label{eq:rate}
\end{align}
\end{proposition}


\end{document}